 \documentclass[11pt, reqno]{amsart}
\usepackage{amsmath,amsthm,amssymb, graphicx}
\theoremstyle{plain}
\newtheorem{thm}{Theorem}[section]
\newtheorem{prop}[thm]{Proposition}
\newtheorem{lem}[thm]{Lemma}
\newtheorem{cor}[thm]{Corollary}

\theoremstyle{definition}

\newtheorem{rem}[thm]{Remark}
\newtheorem{defn}[thm]{Definition}
\newtheorem{eg}[thm]{Example}
\newtheorem{subtitle}[thm]{}
\newtheorem{ex}{Exercise}[section]
\numberwithin{equation}{section}

\def\a{\alpha}

\def\d{\delta}
\def\D{\triangle}

\def\K{\nabla}

\def\n{\,\vert\,}
\def\o{\theta}

\def\co{{\mathcal{O}}}

\def\li{\langle}
\def\ri{\rangle}
\def\n{\ \vert\ }

\def\bs{\bigskip}
\def\ms{\medskip}
\def\ss{\smallskip}

\def\ni{\noindent}
\def\ti{\tilde}
\def\p{\partial}

\def\II{{\rm II\/}}

\def\R{\mathbb{R} }
\def\C{\mathbb{C}}

\newcommand{\beg}{\begin{eg}}
\newcommand{\eeg}{\end{eg}}
\newcommand{\bthm}{\begin{thm}}
\newcommand{\ethm}{\end{thm}}
\newcommand{\bprop}{\begin{prop}}
\newcommand{\eprop}{\end{prop}}
\newcommand{\bcor}{\begin{cor}}
\newcommand{\ecor}{\end{cor}}
\newcommand{\blem}{\begin{lem}}
\newcommand{\elem}{\end{lem}}
\newcommand{\bca}{\begin{cases}}
\newcommand{\eca}{\end{cases}}
\newcommand{\brem}{\begin{rem}}
\newcommand{\erem}{\end{rem}}
\newcommand{\bpm}{\begin{pmatrix}}
\newcommand{\epm}{\end{pmatrix}}
\newcommand{\bbm}{\begin{bmatrix}}
\newcommand{\ebm}{\end{bmatrix}}
\newcommand{\bvm}{\begin{vmatrix}}
\newcommand{\evm}{\end{vmatrix}}
\newcommand{\bdefn}{\begin{defn}}
\newcommand{\edefn}{\end{defn}}
\newcommand{\bsub}{\begin{subtitle}}
\newcommand{\esub}{\end{subtitle}}
\newcommand{\bex}{\begin{ex}}
\newcommand{\eex}{\end{ex}}
\newcommand{\ben}{\begin{enumerate}}
\newcommand{\een}{\end{enumerate}}

\newcommand{\balign}{\begin{align}}
\newcommand{\ealign}{\end{align}}
\newcommand{\baligns}{\begin{align*}}
\newcommand{\ealigns}{\end{align*}}
\newcommand{\beq}{\begin{equation}}
\newcommand{\eeq}{\end{equation}}

\newcommand{\vn}{{\bf n}}
\newcommand{\Id}{{\rm Id\/}}

\begin{document}

\title[Mean curvature flow]
{The mean curvature flow for isoparametric submanifolds}
\author{Xiaobo Liu$^\ast$}\thanks{$^\ast$Research was partially supported by NSF grant DMS-0505835}
\address{Department of Mathematics\\ University of Notre Dame, Notre Dame, IN 46566}
\email{xliu3@nd.edu}
\author{Chuu-Lian Terng$^\dag$}\thanks{$^\dag$Research supported
in  part by NSF Grant DMS-0529756}
\address{Department of Mathematics\\
University of California at Irvine, Irvine, CA 92697-3875}
\email{cterng@math.uci.edu}
%\ms \hskip 3in \today

\begin{abstract}

A submanifold in space forms is {\it isoparametric\/} if the normal bundle is flat and principal curvatures along any parallel normal fields are constant.
We study the mean curvature flow with initial data an isoparametric submanifold in Euclidean space and sphere. We show that the mean curvature flow preserves the isoparametric condition, develops singularities in finite time, and
converges in finite time to a smooth submanifold of lower dimension. We also give a precise description of the collapsing.

\end{abstract}

\maketitle

\section{ Introduction}
\ms

The mean curvature flow (abbreviated as MCF) of a submanifold $M \subset \mathbb{R}^{N}$
over an interval $I$ is a map
$ f: I \times M \longrightarrow \mathbb{R}^{N} $
such that for all $t \in I$ and $x \in M$, $\frac{\partial}{\partial t} f(t, x)$ is equal
to the mean curvature vector  of $M(t)= f(t, M)$ at the point
$x(t) = f(t, x)$.
Mean curvature flows of convex hypersurfaces have been extensively studied in the literature
(cf. \cite{GH}, \cite{Hu}).
An exposition of the work in this area was given in the book \cite{Zhu}.
Comparatively, the behavior of mean curvature flows of submanifolds with higher codimension
is less understood (cf. \cite{W}). This is partly due to the lack of understanding of collapsing
and the formation of singularities of the flow equations in the higher codimensional case.

A submanifold $M$ of a Riemannian manifold is {\it isoparametric} if
its normal bundle is flat and principal curvatures along any parallel normal
vector field are constant. The codimension of $M$ is called the {\it rank} of $M$.
 An isoparametric submanifold $M$ in $\R^N$ is {\it full} if it is not contained in any proper
 hyperplane, and
is {\it irreducible\/} if it is not a product of two isoparametric submanifolds.
We refer to \cite{Terng} for the basic properties and structure
theories for isoparametric submanifolds.
Principal orbits of isotropy representations of symmetric spaces are isoparametric,  they are the only compact homogeneous isoparametric submanifolds in Euclidean space
(cf. \cite{PT}), and are called {\it generalized flag manifolds}.
There are also infinite families of non-homogeneous isoparametric submanifolds which arise from
representations of Clifford algebras (cf. \cite{FKM}).
All these non-homogeneous examples have rank 2. A theorem of Thorbergsson \cite{Th} asserts
that compact full irreducible isoparametric submanifolds with rank bigger than 2 are always
homogeneous.

A complete isoparametric submanifold of $\mathbb{R}^{N}$ can be
decomposed as the product of a compact, irreducible, isoparametric
submanifold and a subspace of $\mathbb{R}^{N}$. Since mean
curvature flows with affine subspaces of $\mathbb{R}^{N}$ as
initial data is trivial and the mean curvature flow starting from
a product submanifold stays as product, we will only consider
compact, full, irreducible isoparametric submanifolds.

Let $M$ be an isoparametric submanifold of $\R^N$, and $\xi$ a parallel normal vector field
on $M$. Then
$M_\xi=\{x+ \xi(x)\n x\in M\}$
is again a smooth submanifold (may have higher codimension), and the map $M\to M_\xi$ is either a
diffeomorphism or a fibration with a generalized flag manifold as fiber.
The family of these parallel sets forms a singular foliation of the ambient Euclidean space $\R^N$.
Top dimensional leaves are all isoparametric in $\mathbb{R}^{N}$, and they are called {\it parallel
isoparametric submanifolds}.
Lower dimensional leaves are no longer isoparametric, and they are called
{\it focal submanifolds\/}.

We show that if $f:M\times [0, T)\to \R^N$ is a solution of the MCF in $\R^N$ with
$f(\cdot, 0)$ isoparametric then $f(\cdot, t)$ is isoparametric for all $t\in [0,T)$, i.e.,
the MCF preserves isoparametric condition.  This reduces the MCF to a system of ordinary
differential equations. There is a Weyl group $W$ associated to each isoparametric submanifold $M$ that acts on the normal plane $p+ \nu_p M$. The ODE given by the mean curvature flow with initial data an isoparametric submanifold is given by a vector field $H$ smoothly defined on the interior of the Weyl chamber $C$ of $W$ but blows up at the boundary of $C$.  However, we can use generators of $W$-invariant polynomials to change coordinate so that the vector field $H$ becomes a polynomial vector field and its flows can be solved explicitly.

Every compact isoparametric submanifold is contained in a sphere. This sphere is
also foliated by parallel isoparametric submanifolds and focal submanifolds.
Each isoparametric foliation contains a unique isoparametric submanifold which is a minimal
submanifold of this sphere. The mean curvature flow in $\R^N$ with initial data a minimal submanifold
in $S^{N-1}$ behaves like the mean curvature flow of  a sphere, i.e. it just shrinks homothetically
along the radial direction and collapses to a point in finite time
(cf. Lemma \ref{miniso}). If $M$ is an isoparametric submanifold in $\R^N$
which is not minimal in the sphere, then its mean curvature flow will converge to a focal
submanifold $F$ of positive dimension (cf. Corollary \ref{cor:miniso}).
In fact, $M$ is a fibration over $F$ with each fiber a homogeneous isoparametric submanifold of a
lower dimensional Euclidean space.
Each fiber of this fibration collapses
to a point under the mean curvature flow in a finite time.

We summarize some of the main results of this paper in the following Theorem (cf. Theorem \ref{thm:finiteconv}, Theorem \ref{thm:focalconv}, and Proposition \ref{bm}):

\begin{thm}\label{bo}
The mean curvature flow in $\R^N$ with initial data a compact isoparametric submanifold
\ben
\item converges to a focal submanifold in finite time $T$,
\item if the fibration from the initial isoparametric submanifold to the limiting focal
submanifold is a sphere fibration (this is the generic case), then the mean curvature flow
$M(t)$ has type I singularity, i.e., there is a constant $c_0$ such that
$||\II(t)||_\infty^2(T-t)\leq c_0$ for all $t\in [0,T)$, where $||\II(t)||_\infty$ is
the sup norm of the second fundamental form of $M(t)$,
\item every focal submanifold is the limit of
the mean curvature flow with some parallel isoparametric submanifold as initial data,
\item  if $M_1$ and $M_2$ are distinct parallel full isoparametric submanifolds in
$\R^N$ that lie in the same sphere. Then the mean curvature flows in $\R^N$ with initial data
$M_1$ and $M_2$ collapse to two distinct focal submanifolds.
\een
\end{thm}

 The mean curvature flow in $S^{N-1}$ with initial data an isoparametric  submanifold behaves very similarly to
the Euclidean mean curvature flow. In particular we have the following theorem:

\begin{thm} \label{thm:MCFisopS}
Let $M$ be an isoparametric submanifold of $S^{N-1}$.  Then the mean curvature flow in
$S^{N-1}$ with $M$ as initial data
\ben
\item  is constant if $M$ is minimal in $S^{N-1}$, or
\item converges to a focal submanifold of positive dimension in finite
time if $M$ is not minimal.
\een
\end{thm}

An isometric action of $G$ on a Riemannian manifold $N$ is {\it polar\/} if there exists a closed embedded submanifold $\Sigma$ of $N$ that meets all $G$-orbits and meets orthogonally. Such $\Sigma$ is called a {\it section\/} of the polar action.  Principal orbits of a polar action in $\R^n$ and $S^n$ are isoparametric (cf. \cite{PT}).  We prove that if the $G$-action on $N$ is polar then the mean curvature flow preserves $G$-orbits and the flow becomes an ordinary differential equation on the section $\Sigma$.  We expect that methods developed in this paper can be applied to study mean curvature flows for orbits of polar actions with flat sections in symmetric spaces.

This paper is organized as follows:  We give a brief review of properties of isoparametric submanifolds that are needed in section \ref{sec:prelim}, present proofs of results stated in Theorem \ref{bo} in section \ref{sec:Basic}, construct explicit solutions of the MCF in $\R^N$ with initial data an isoparametric submanifold in section \ref{sec:solution}.  Since focal submanifolds are smooth manifolds, we can consider their mean curvature
flow. Most properties of the mean curvature flows for isoparametric submanifolds
also hold for focal submanifolds. This will be briefly discussed in section \ref{sec:focal}.
We describe MCF in spheres with initial data an isoparametric submanifold in spheres in section \ref{sec:sphere}, and in the last section we discuss MCF in  a Riemannian manifold $N$ with initial data a principal orbit of a polar action on $N$.

The authors like to thank Mu-Tao Wang, Yng-Ing Lee, and Mao-Pei Tsui for discussions on types of singularities of the mean curvature flows.

\bs
\section{Preliminaries}
\label{sec:prelim}

Geometric and topological properties of isoparametric submanifolds can be found in
\cite{Terng}. In this section we
briefly review the properties which will be used in this paper.
Let $M \subset \mathbb{R}^{N}$ be a full compact
isoparametric submanifold of rank $k$.

\subsection{Curvature spheres and curvature normals}

\hspace{100pt}\\
The tangent bundle of $M$ can be decomposed into
orthogonal sums of {\it curvature distributions} $\{ E_{i} \mid i=1, \cdots, g\}$
for some integer $g>0$. At each point of $M$, $E_{i}$ is a common eigenspace
of the shape operators of $M$ at that point.
 There are parallel normal vector fields $\vn_{i}$
such that the shape operator $A_{\xi}$ has the property
\[ A_{\xi} |_{E_{i}} = \li \xi, \vn_{i}\ri  \Id_{E_i} \]
for all normal vector $\xi$.
Each vector field $\vn_{i}$ is called the {\it curvature normal} of $E_{i}$. The
rank of $E_{i}$ is called the {\it multiplicity} of $\vn_{i}$,
which will be  denoted by $m_{i}$.
Each $E_{i}$ is an integrable distribution whose leaves are $m_{i}$-dimensional
round spheres with radius $1/\|\vn_{i}\|$. Such spheres are called
{\it curvature spheres}.

\subsection{Weyl group associated to $M$}\label{ae}

\hspace{100pt}\\
For each $i \in \{1, \cdots, g\}$, let $\sigma_{i}(x)$ be the antipodal point in the
$i$-th curvature sphere passing through $x$. Then $\sigma_{i}$ is an involution on $M$.
The group $W$ generated by $\sigma_{1}, \cdots, \sigma_{g}$ is  a crystallographic {\it Coxeter group}.
It is known that $M$ is irreducible if and only if $W$ is irreducible. For each $x \in M$,
$W$ also acts as a reflection group
on the affine normal space $x+\nu_{x}M$
generated by reflections along hyperplanes
\[ L_{i} := \{ x+\xi \mid \xi \in \nu_{x}M, \,\,\, 1- \li \xi, \vn_{i}(x)\ri = 0 \} \]
for $i=1, \cdots, g$.

The intersection $\bigcap_{i=1}^{g} L_{i}$ consists of a single constant point
which is denoted by  $a$. Then $M$ is contained in a sphere which is centered at $a$.
Without loss of generality, we always {\bf assume} that $a=0$, i.e.,  $M$ is contained in a sphere
centered at the origin of $\mathbb{R}^{N}$.  This condition is equivalent to
\begin{equation} \label{eqn:centerO}
 \li -x, \, \vn_{i}(x)\ri \,\, =  \,\, 1
\end{equation}
for all $x \in M$ and $i=1, \cdots, g$ (cf. \cite[Corrolary 1.17]{Terng}).

%\ss
%\ni{\bf 1.3. Parallel submanifold}\label{af}

\subsection{Parallel submanifold}

\hspace{100pt} \\
 For any parallel normal vector field $\xi$ on $M$, define
\[ M_{\xi}:= \{ x+\xi(x) \mid x \in M \}. \]
If
\begin{equation} \label{eqn:noncollcond}
 1- \li\xi(x), \vn_{i}(x)\ri \neq 0
\end{equation}
 for $i=1, \cdots, g$, then
$M_{\xi}$ is again an isoparametric submanifold with the same dimension as $M$.
$M_{\xi}$ is called the {\it parallel isoparametric submanifold} of $M$ defined by $\xi$.
The curvature normals of $M_{\xi}$ at the point $x+ \xi(x)$ are given
by
\[ \frac{\vn_{i}(x) }{1- \li\xi(x), \vn_{i}(x)\ri} \]
with same multiplicities $m_{i}$ for $i=1, \cdots, g$.
The mean curvature vector of $M_{\xi}$ at $x+ \xi(x)$ is given
by
\begin{equation} \label{eqn:MCV}
 H(x+\xi(x)) = \sum_{i=1}^{g} \frac{m_{i}\vn_{i}(x) }{1- \li\xi(x), \vn_{i}(x)\ri}.
\end{equation}

When condition \eqref{eqn:noncollcond} fails,  $M_{\xi}$ is still
a smooth submanifold of $\mathbb{R}^{N}$, but it is no longer
isoparametric. This submanifold is called
a focal submanifold of $M$. The dimension of $M_{\xi}$ is strictly
smaller than that of $M$. The map
\[ \begin{array}{crcl}
\pi:& M & \longrightarrow & M_{\xi} \\
    & x & \mapsto & x+\xi(x)
\end{array} \]
is a fibration over $M_{\xi}$ with each fiber an isoparametric submanifold in the normal space
of $M_{\xi}$ at $\pi(x)$.
 In fact, fix $x_0\in M$, let $C$ denote the Weyl chamber of $W$ on $x_0+ \nu_{x_0}M$
 containing $x_0$, i.e.,
 $$C=\{x_0+\xi\n \xi\in \nu_{x_0}M, \, \li \xi, n_i\ri <1\}.$$
If $y_0=x_0+\xi(x_0)$ lies in the boundary of $C$ and $y_0\not=0$, then the fiber of the
fibration $M\to M_\xi$ is a generalized flag manifold with Weyl group $W_{y_0}$,
the isotropy subgroup of $W$ at $y_0$.

 For any parallel normal vector field $\xi$ on $M$, the intersection of $M_{\xi}$ with
$x + \nu_{x}M$ is an orbit of $W$. In particular,  if $M_{\xi}$
is a parallel isoparametric submanifold, then it intersects each open Weyl chamber of
$W$ exactly once. Moreover, $M_{\xi}$ is a focal submanifold if and only if
$ x+\xi(x)$ is contained in $\bigcup_{i=1}^{g} L_{i}$.

%\ss\ni {\bf 1.4. Isoparametric map and $W$-invariant polynomials}\label{ag}
%\ss

\subsection{Isoparametric map and $W$-invariant polynomials}

\hspace{100pt} \\
Given a $W$-invariant polynomial $f$ on $V=x_0+\nu_{x_0}M$, there is a unique extension
$\Psi(f)$ on $\R^N$ such that $\Psi(f)$ is constant along any parallel submanifold $M_\xi$
and $\Psi(f)\n_{V}= f$.  Moreover, $\Psi(f)$ is also a polynomial.

Let $\ti \D$ and $\D$ denote the Laplacian in $\R^N$ and $V$ respectively.
Then by Lemma 3.2 of \cite{Terng},
\beq\label{eqn:mcfpoly}
F(x)=\ti \D \Psi (f)(x)- \Psi(\D f)(x)= \sum_{i=1}^g \frac{m_i \li \K f(x), \vn_i\ri}{\li x, \vn_i\ri}.
\eeq
is a polynomial on $\R^N$ and is constant along parallel submanifolds of $M$.
Moreover,  if $f$ is a homogeneous $W$-invariant polynomial of degree $m$ on $V$,
then  $F$ is a homogeneous polynomial of degree $m-2$ on $\R^N$.

\bs
\section{Mean curvature flows for general isoparametric submanifolds}
\label{sec:Basic}

Let $M$ be an isoparametric submanifold of $\R^N$, fix $x_0\in M$, and $W$ the Coxeter group associated to $M$.
We prove that the MCF stays isoparametric and the MCF equation becomes a flow equation of a vector field $H$ defined in the interior of the Weyl chamber of $W$ containing $x_0$ in $x_0+ \nu_{x_0}M$ and the vector field $H$ tends to infinity at the boundary of the Weyl chamber.  We prove that solutions of the ODE $x'= H(x)$ only exists for finite time. To see the finer structure of the behavior of the blow-up of MCF,  we use $W$-invariant polynomials to construct a new coordinate system for the Weyl chamber so that the corresponding vector field $H$ becomes a polynomial vector field. We then analyse the behavior of flows of this polynomial vector field to obtain informations on the collapsing of the MCF.

Fix $x_{0} \in M$.  Since $\nu M$
is globally flat, we can identify a vector $v\in \nu_{x_0}M$ and the unique parallel normal field $\hat v$
along $M$ defined by $\hat v(x_0)=v$.
Let $\vn_{i}$ be curvature normals of $M$ with multiplicity $m_{i}$ for $i=1, \ldots, g$.
We may view $\vn_{i}$ either as a global parallel normal vector field along $M$ or
an element in $\nu_{x_0}M$. The precise meaning should be clear from the context.

Let $\xi(t) \in \nu_{x_0}M$ be a one parameter family of normal vectors satisfying the flow equation
\begin{equation} \label{eqn:mcf}
 \dot{\xi}(t) = \sum_{i=1}^{g} \frac{m_{i}\vn_i}{1- \li\xi(t), \vn_{i}\ri}, \qquad
    \xi(0) = 0.
\end{equation}
It follows from  \eqref{eqn:MCV} that $\xi$ is a solution of \eqref{eqn:mcf} if and only if
 the one parameter family of parallel submanifolds $M(t) := M_{\xi(t)}$ satisfy the
{\it mean curvature flow} equation  with $M(0) = M$.
In other words, the MCF preserves the isoparametric condition:

\bprop
If $f:M\times [0,T)\to \R^N$ satisfies the mean curvature flow in $\R^N$ and $f(\cdot, 0)$ is isoparametric,
then $f(\cdot, t)$ is isoparametric for all $t\in [0,T)$.
\eprop

Equation \eqref{eqn:mcf} does not make sense if $\li\xi(t), \vn_{i}\ri = 1$ for some $i$.
We will only study the flow equation under the condition
\[ \li\xi(t), \vn_{i}\ri \,\, <  \,\, 1 \]
for all $i=1, \ldots, g$. In other words, we require that $x_{0} + \xi(t)$ stays in the same
Weyl chamber as $x_{0}$ for all $t$. Under this condition, all $M(t)$ are still isoparametric.

Note that \eqref{eqn:mcf} is a system of non-linear ODE given by a vector field  defined on the
Weyl chamber $C$ containing $x_0$ and the vector field blows up along the boundary of $C$.
The study of MCF with isoparametric submanifolds as initial data reduces to the study of
solutions of this ODE system.

\begin{thm} \label{thm:radialMCV}
Let $M\subset S^{N-1}(r_0)$ be an $n$-dimensional isoparametric submanifold in $\R^N$, and $x_0\in M$.
If $\xi(t)$ satisfies the mean curvature flow equation \eqref{eqn:mcf}, then $x(t)=x_0+\xi(t)$
satisfies
\beq\label{aa}
x'(t)= -\sum_{i=1}^g \frac{m_i \vn_i}{\li x(t), \vn_i\ri},
\eeq
with $x(0)=x_0$.
Let $H(t)$ be the mean curvature vector of $M(t)=M_{\xi(t)}$ at the point $x(t)$. Then
\ben
\item[(a)] $\li x(t), \, H(t)\ri = -  n$,
\item[(b)] $|| x(t) ||^{2} = || x(0)||^{2} - 2nt$.
\een
\end{thm}

\begin{proof}
By equation \eqref{eqn:centerO},
$$ \li x(t), \, \vn_{i}\ri  =\li x(0), \vn_{i}\ri  + \li \xi(t), \vn_{i}\ri \nonumber = -1+ \li \xi(t), \vn_{i}\ri$$
for all $i=1, \cdots, g$.
Since
\begin{equation} \label{eqn:MCVcenterO}
H(t) = - \sum_{i=1}^{g} \frac{m_{i}\vn_i}{\li x(t), \vn_{i}\ri},
\end{equation}
we have $\li x(t), \, H(t)\ri = -  \sum_{i=1}^{g} m_{i} = -n$.
This proves part (a). Part (b) follows from integrating the following formula
\[ \frac{d}{dt} \|x(t)\|^{2}
\,\,  = \,\, 2 \li x(t), \, x^{\prime}(t)\ri
\,\,=\,\, 2 \li x(t), \, H(t)\ri
\,\,= \,\, -  2n.\]
\end{proof}

Hence we have
\begin{cor} \label{cor:finiteInt}
The  mean curvature flow in $\R^N$ with initial data an isoparametric submanifold in $S^{N-1}(r_0)$ exists only for finite time with maximal interval $[0, T)$, where
$0<T \leq T_{0}=\frac{r_0^2}{ 2n}$.
\end{cor}

The following Theorem is the key in proving
\ben
\item the limits of two flows of \eqref{aa} have two different limit on the boundary $\p C$ of the Weyl chamber $C$,
\item every point of $\p C$ is a limit of some flow of \eqref{aa}.
\een

\bthm\label{ba}
Let $M$ be a compact isoparametric submanifold in $\R^N$, $W$ the Weyl group associated to $M$, $x_0\in M$ a fixed point, and $V=x_0+\nu_{x_0}M$.
Let $P_1, \ldots, P_k$ be a set of generators of the ring $\R[V]^W$ of $W$-invariant polynomials
on $V$ such that $P_i$ are homogeneous polynomials of degree $s_i$ with
$P_1(x)=||x||^2$ and $s_1\leq s_2\cdots \leq s_k$, and $C$ the Weyl chamber of $W$
containing $x_0$ in $V$.  Let $P:\bar C\to \R^k$ be the map defined by
$P(x)= (P_1(x), \ldots, P_k(x))$.  Then
 $P$ is a homeomorphism from $\bar C$ to a closed subset $B=P(\bar C)$. Moreover, there
 is a polynomial map
 $$\eta=(\eta_1, \ldots, \eta_k):\R^k\to \R^k$$
 with $\eta_1=-2n$ and $\eta_j$ is a polynomial in $\eta_1, \ldots, \eta_{j-1}$ such that
 if $x:[0,T)\to C$ is a solution of \eqref{aa}, then $y(t)= P(x(t))$ is a solution of
$$y'(t)= \eta(y(t))= (\eta_1(y(t)), \ldots, \eta_k(y(t))).$$
\ethm

\begin{proof}
Since each orbit of $W$ intersect $\overline{C}$
exactly once and the algebra of invariant polynomials separate $W$ orbits,
the map
\[ \begin{array}{crcl}
    P: & \overline{C}\cap D(1) & \longrightarrow & \mathbb{R}^{k} \\
        & x & \longmapsto &  (P_{1}(x), \cdots, P_{k}(x))
        \end{array}   \]
is an injective continuous  map.
Since $P$ is injective and proper,  $B=P(\overline{C})$ is a closed subset of $\mathbb{R}^{k}$ and
 $P$ is a homeomorphism from $\overline{C}$ to $B$.

The Coxeter group $W$ acts on $V$. If $f$ is a $W$-invariant
homogeneous polynomial of degree $j$ on $V$,
then by equation \eqref{eqn:mcfpoly},
$$ F(x) := \sum_{i=1}^{g} m_{i}\,\, \frac{\li \nabla f(x), \, \vn_{i}\ri}{\li x, \, \vn_{i}\ri }$$
is a $W$-invariant homogenous polynomial of degree $j-2$.
Let $x(t)$ be a solution of \eqref{aa}, and $f(t)=f(x(t))$.
By equation \eqref{eqn:MCVcenterO},
\[ f^{\prime}(t) = \li \nabla f(x(t)), x^{\prime}(t)\ri
    = \li \nabla f(x(t)), H(t)\ri = - F(x(t)). \]
Hence $f^{\prime}(t)$ is the value of a $W$-invariant polynomial
of degree $k-2$ evaluated at $x(t)$.
In particular, $\frac{d^{j}}{dt^{j}} f(t) = 0$ if $j>k/2$.
Therefore $f(t)$ is a polynomial in $t$.

Let $y(t)= (y_1(t), \ldots, y_k(t))= P(x(t))$, and
 $$F_i(x)= \sum_{i=1}^g m_{i} \, \frac{\li \K P_i(x), \vn_i\ri}{\li x, \vn_i\ri}.$$
By equation \eqref{eqn:mcfpoly}, $F_i$ is a $W$-invariant homogeneous polynomial on
$V$ of degree $s_i-2$. Since $\R[V]^W=\R[P_1, \ldots, P_k]$,
    $$F_i= - \eta_i(P_1, \ldots, P_{i-1})$$
for some polynomial $\eta_i$.  But we have shown above that $y_{i}'(t)= - F_i(x(t))$, so
%   $$y_i'(t)= - F_{i}(x(t))= f_i(P_1(x(t)), \ldots, P_{i-1}(x(t)))
%   = \eta_i(y_1(t), \ldots, y_{i-1}(t)).$$
$$y_i'(t)= - F_{i}(x(t)) = \eta_i(y_1(t), \ldots, y_{i-1}(t)).$$

This shows that $y(t)$ is an integral curve of the polynomial vector field $\eta$ on $\R^k$.
Since $y_1(t)= ||x(0)||^2 -2nt$, solution $y$ can be solved explicitly by integrations.
 \end{proof}

The MCF equation \eqref{aa} is given by the vector field
 $$H(x)= -\sum_{i=1}^g \frac{m_i \vn_i}{\li x, \vn_i\ri},$$ which is smoothly defined
 on the Weyl chamber $C$ of $x_0+\nu_{x_0}M$ and blows up at the boundary $\p C$.
 If we use generators of $W$-invariant polynomials on $x_0+\nu_{x_0}M$ to change coordinates
 to $P$ as in Theorem \ref{ba}, then the vector field $H$ becomes the polynomial vector field $\eta$ on $P(C)$.
 Moreover, the flow of $\eta$  can be solved explicitly and globally. Then apply $P^{-1}$
 to flows of $\eta$ to get flows of \eqref{aa}.

\begin{thm} \label{thm:finiteconv}
For any compact isoparametric submanifold $M$ in $\R^N$,  the mean curvature flow always converges
to a focal submanifold at a finite time.  Moreover, if $M_{1}$ and $M_{2}$ are parallel full
isoparametric submanifolds which
 are contained in the same sphere, then mean curvature flows with initial data $M_1$ and $M_2$
 never intersect and they
 converge to two distinct focal submanifolds.
\end{thm}

\begin{proof}
 We use the same notation as Theorem \ref{ba}.
Let $x:[0, T)\to C$ be the maximal interval for a solution of the
mean curvature flow equation \eqref{aa}. Note that
$P_i(t)=P_{i}(x(t))$ are well defined since $P_{i}(t)$ are
polynomials in $t$. Therefore the mean curvature flow of $x_{0}
\in M$ must converge to $P^{-1}(P_{1}(T), \cdots, P_{k}(T))$ which
lies on the boundary of $\overline{C}$. The mean curvature flow of
$M$ then converges to the focal submanifold passing through this
point.

We may assume that $x_i(0)$ lies in the unit sphere.  Let $T_i$ denote the maximum time for the solution $x_i(t)$. If $T_1\not= T_2$, then $||x_i(t)||^2= 1- 2n t$, so $\lim_{t\to T_1^-}||x_1(t)||^2\not=\lim_{t\to T_2^-} ||x_2(t)||^2$.  If $T_1=T_2= T$, then since $x_i(t)$ are solutions of \eqref{aa} and $\li x_i(t), n_j\ri <0$, we have
\begin{equation}\label{bp}
\frac{1}{2} \frac{d}{dt}\, \big| ||x_1(t)-x_2(t)||^2= \sum_{i=1}^g m_i\, \frac{\li x_1(t)-x_2(t), \vn_i\ri^2}{\li x_1(t), n_i\ri \li x_2(t), \vn_i\ri} \geq 0.
\end{equation}
This implies that $||x_1(t)-x_2(t)||^2$ increases in $t\in [0, T)$, hence
$$\lim_{t\to T^-} x_1(t)\not= \lim_{t\to T^-} x_2(t).$$
\end{proof}

\begin{thm}\label{thm:focalconv}
Every focal submanifold is a limit of the mean curvature flow of certain isoparametric submanifold.
\end{thm}

We need a couple Lemmas.
First a simple Lemma on scaling and the proof is obvious:

\blem\label{ac}
If $f:M\times [0,T)\to \R^N$ is a solution to the mean curvature flow with $f(x,0)= f_0(x)$,
then given any $r\not=0$, $\ti f(x,t)= rf(x, r^{-2}t)$ is a solution with $\ti f(x,0)= rf_0(x)$.
\elem

\begin{lem} \label{miniso}
Let $f:M^{n} \longrightarrow S^{N-1}(r_{0})$ be an immersed minimal submanifold of
a sphere with radius $r_{0}$.
For any $x \in M$, the solution to the mean curvature flow equation  in $\mathbb{R}^{N}$
with initial data $M$ is given by
$$F(x, t) = \sqrt{1 - (2nt/r_{0}^{2}) } \, f(x).$$
In particularly, the mean curvature flow of $M$  shrinks to a point homothetically
 in finite time
$T_{0} \,\, = \,\, r_{0}^{2}/(2n)$.
\end{lem}

\begin{proof}
For minimal submanifolds of the sphere $S^{N-1}(r)$ with radius $r$,
the mean curvature vector at a point $x$ is $- n x/r^{2}$.
Let $F(x, t)= r(t) f(x)$ for $x \in M$ with $r(t) \geq 0$. Then the mean curvature
vector field of $F(\cdot, t)$ at point $x$ is given by
$- \frac{n}{ r_{0}^{2} r(t)} \, f(x)$.
So $F(x, t)$ satisfies the mean curvature flow equation for $f$ if and only if
\[ r^{\prime}(t) = - n/( r_{0}^{2} r(t)) \hspace{20pt}
    {\rm and} \hspace{20pt} r(0)=1.\]
It follows that $ r(t) = \sqrt{1 - (2nt/r_{0}^{2})}$.
\end{proof}

\ms
\ni {\bf Proof of Theorem \ref{thm:focalconv}}

In each isoparametric family, there exists a unique isoparametric submanifold $M\subset S^{N-1}(1)$,
which is minimal in $S^{N-1}(1)$. Let $x_0\in M$. The mean curvature flow for minimal
 submanifold in spheres can be solved explicitly as in Lemma \ref{miniso}, i.e.,
 $x(t)= \sqrt{1-2nt}\, x_0$ is a solution of \eqref{aa} and $x(t)\in C$ for all
 $t\in [0,\frac{1}{2n})$.

 Recall that integral curves of $H(x)= -\sum_{i=1}^g \frac{m_i\vn_i}{\li x, \vn_i\ri}$ map to
 integral curves of the polynomial vector field $\eta$ under the homeomorphism $P$ defined in
 Theorem \ref{ba}.  Since the integral curve starting from $x_0$ lies in $C$, the flow of $\eta$
 starting at $P(x_0)$ lies in $P(C)$.  But $-\eta$ is a polynomial vector field and the one-parameter
 subgroup $\phi_t$ generated by $-\eta$ is a globally defined polynomial map.  So there exists $\d>0$
 and an open subset $\mathcal U$ of $P(\bar C)$ such that $\mathcal U$ contains the origin and
 $\phi_t(z)\in P(C)$ for $t\in (0,\d)$ and
 $z\in {\mathcal U}$.  This shows that the flow of $-\eta$ starting at the boundary of $\mathcal U$
 points inward in $P(C)$.   Hence any boundary point of $P^{-1}({\mathcal U})$ is a limit of some
 MCF with initial data in $C$. It follows from Lemma \ref{ac} that any focal submanifold can be a
 limit of some MCF with some initial isoparametric submanifold.
 $\Box$

As consequence of Theorem  \ref{thm:finiteconv} and Lemma \ref{miniso}, we have

\begin{cor} \label{cor:miniso}
Let $M$ be an isoparametric submanifold.
The mean curvature flow of $M$ converges to a point if and only
if it is minimal in the sphere containing it.
\end{cor}

Below we describe the rate of collapsing of the MCF for isoparametric submanifolds.  Recall that a MCF, $M_t$, collapses at time $T<\infty$ is said to have {\it type I singularity\/} (cf. \cite{W}) if there is a constant $c_0$ such that
$$||\II(t)||_\infty^2(T-t)\leq c_0$$
for all $t\in [0,T)$, where $||\II(t)||_\infty$ is the sup norm of the second fundamental form for $M_t$.

\bprop\label{bm}
Let  $x(t)$ be a solution of the MCF \eqref{aa}, and $T$ is the maxmial time.  Then
\ben
\item $x(T):=\lim_{t\to T^-} x(t)$ exists and belong to the boundary $\p C$ of the Weyl chamber $C$,
\item $\lim_{t\to T^-} \frac{||x(t)-x(T)||^2}{T-t} = 2m$, where $m=\dim(M_{x(0)})-\dim(M_{x(T)})$,
\item if $x(T)$ lies in a highest dimensional stratum of $\p C$, then the MCF has type I singularity.
\een
\eprop

\begin{proof}
We have proved (1) in Theorem \ref{thm:finiteconv}. Statement (2) follows from the L'Hopital law:
\begin{align*}
& \lim_{t\to T^-} \frac{||x(t)-x(T)||^2}{T-t}
=\lim_{t\to T^{-}} \frac{2 \li x(t)-x(T), \,  x'(t)\ri}{-1} \\
&= 2\lim_{t\to T^{-1}} \sum_{i=1}^{g} \li x(t)-x(T), \, \frac{m_i \vn_i}{\li x(t), \, \vn_i\ri}\ri\\
&=2 \lim_{t\to T^{-1}} \sum_{i\not\in I} m_i \frac{\li x(t)-x(T), \, \vn_i\ri}{\li x(t), \, \vn_i\ri}
 + 2\sum_{i\in I} m_i\frac{\li x(t), \, \vn_i\ri }{\li x(t), \, \vn_i\ri} = 2 \sum_{i\in I} m_i,
\end{align*}
which is the dimension of the fiber of $M_{x(0)}\to M_{x(T)}$.
Here $I=\{1\leq i\leq g\n \li x(T), \, \vn_i\ri=0\}$.

We now prove statement (3).
If $x(T)$ lies in a highest dimensional stratum of $\p C$, then there exists a  unique $i$ such that $x(T)$ lies
in the hyperplane defined by $\vn_i$, i.e.,  $\li x(T), \vn_i\ri=0$.  We may assume $i=1$.
Note that the norm square of the second fundamental form of $M_{x(t)}$ satisfies
\begin{align*}
& ||\II(x(t))||^2(T-t) \leq \sum_{i=1}^g \frac{m_i ||\vn_i||^2}{\li x(t), \vn_i\ri^2}(T-t) \\
&= \frac{m_1||\vn_1||^2(T-t)}{\li x(t), \vn_1\ri^2}+\sum_{i=2}^g \frac{m_i||\vn_i||^2}{\li x(t), \vn_i\ri^2}(T-t).
\end{align*}
As $t\to T^-$, the second term tends to zero because $\li x(T), \vn_i\ri\not=0$ for
all $i\geq 2$, and by the l'Hopital law the first term has the same limit as
$$\frac{-m_1||\vn_1||^2}{-2\li x(t), \vn_1\ri \sum_{i=1}^g m_i \frac{\li\vn_i, \vn_1\ri}{\li x(t), \vn_i\ri}}.$$
But the denominator tends to $-2m_1||\vn_1||^2$, so the limit is $1/2$.
\end{proof}

We remark that there is an open dense subset $\co$ of the Weyl chamber $C$ such that the solution $x(t)$ of \eqref{aa} with $x(0)\in \co$ converges to a point in a highest dimensional stratum of $\p C$.

\bs
\section{Solutions to the mean curvature flow equation}
\label{sec:solution}

In this section, we use Theorem \ref{ba} to construct explicit solutions of the MCF \eqref{aa} by
selecting a set of generators $P_1, \ldots, P_k$ for the $W$-invariant polynomials and calculating
flows of the polynomial vector field $\eta$.

We use the root system of the Coxeter group given in \cite{GB}.
Let $M$ be a compact, irreducible isoparametric submanifold in $\R^N$, $W$ its Weyl group, and
$\vn_i$ its curvature normals.
Let $\Pi$ denote a set of simple roots of $W$, and $\Delta_{+}$ the set of positive roots
defined by $\Pi$.   Then $\{\R\vn_i\n 1\leq i\leq g\}$ is equal to $\{\R\a\n \a\in \Delta_+\}$.  So
the Weyl chamber $C$ containing $x_0$ is precisely given by
\[ C = \{ x \in V \mid
    -\li x, \a\ri > 0 \,\,\, {\rm for \,\,\, all} \,\,\, \alpha \in \Pi \}. \]
The closure of $C$ is
\[ \overline{C} = \{ x \in V \mid
    -\li x, \a\ri \geq 0 \,\,\, {\rm for \,\,\, all} \,\,\, \alpha \in \Pi \}, \]
and the MCF \eqref{aa} becomes
\begin{equation}\label{ab}
x'(t)=-\sum_{\a\in \Delta^+} \frac{m_\alpha }{\li x(t),\a\ri}\,\, \a
\end{equation}
where $m_{\alpha}$ is the multiplicity of the curvature normal which is parallel to $\a$.
Since \eqref{ab} is invariant under re-scaling of each $\a$, we may
normalize roots of the Coxeter group to be of {\it unit length}.

If $M = M_{1} \times M_{2}$ with $M_{i}$ an isoparametric submanifold of $\mathbb{R}^{N_{i}}$
for $i=1, 2$, then the Weyl group of $M$ is the product of the Weyl groups of $M_1$ and $M_2$
and the mean curvature flow of $M$ is the product
of the mean curvature flows of $M_{1}$ and $M_{2}$.
 So without loss of generality,
we may assume that $M$ is an irreducible isoparametric submanifold.
In the rest of this section, we work out explicit solutions for mean curvature flow equations
for compact isoparametric submanifolds whose Coxeter group are $A_{k}$, $B_{k}$,
$D_{k}$ and $G_{2}$.

\begin{eg} {\bf The $A_{k}$ case}\par

\ss

Suppose that $k \geq 2$.
Let $\{ e_{1}, \cdots e_{k+1} \}$ be the standard orthonormal basis of $\mathbb{R}^{k+1}$
and $(x_{1}, \cdots , x_{k+1})$ the corresponding coordinate.
The set
$\frac{1}{\sqrt{2}} \left(e_{i} - e_{i+1} \right)$ with $1 \leq i \leq k$ is a simple root
system of $A_k$, and
the set of positive roots is
$\frac{1}{\sqrt{2}} \left(e_{i} - e_{j} \right)$ with $1 \leq i<j \leq k+1$.
Let
\[ V:= \{(x_{1}, \cdots, x_{k+1}) \in \mathbb{R}^{k+1} \mid x_{1}+\cdots+x_{k+1} = 0 \}. \]
The normal space of an isoparametric submanifold of type $A_{k}$ can be identified with $V$.
The Coxeter group acts on
$V$ and is generated by all permutations of $\{e_{1}, \cdots, e_{k+1}\}$ .
The open positive Weyl chamber $C$ containing $x_0$ is
\[ C=\{ (x_{1}, \cdots, x_{k+1}) \in V \mid x_{1} < x_{2} < \cdots < x_{k+1}  \}.\]
 Since the multiplicities of
curvature spheres are invariant under the action of the Coxeter group, there is only one
possible multiplicity which we denote by $m$.  So \eqref{ab} is
\beq \label{ak}
x'= - \sum_{i<j} m\, \frac{e_i-e_j}{x_i-x_j},
\eeq
which implies that
\begin{eqnarray} \label{eqn:MCFA}
\frac{1}{m} \frac{d}{dt} x_{i}
&=&  \sum_{j \neq i} \frac{1}{x_{j} - x_{i}}, \quad 1\leq i\leq k+1.
\end{eqnarray}
If $x(0) \in V$, then
$x(t) \in V$ for all $t$. This follows from the simple fact that
$ \frac{d}{dt} (x_{1} + \cdots x_{k+1}) = 0$.

Let $\sigma_{r}$ be the $r$-th elementary symmetric polynomial in $x_{1}, \cdots, x_{k+1}$:
$$\sigma_r=\sum_{1\leq i_1<\cdots < i_r\leq k+1} x_{i_1}\cdots x_{i_r},$$
and $\sigma_{0}=1$.
Let $x(t)$ be a solution of \eqref{ak}, and
$$y_r(t)= \sigma_r(x(t)).$$
We claim that
\beq\label{ak1}
\bca y_2'= n,\\
y_r'= \frac{1}{2}\, m(k-r+3)(k-r+2) y_{r-2}.\eca
\eeq

To see this, we compute
\begin{eqnarray*}
\frac{r!}{m} \frac{d}{dt} y_{r}
&=& \frac{1}{m} \frac{d}{dt} \sum_{i_{1} \neq \cdots \neq i_{r}} x_{i_{1}} \cdots x_{i_{r}} \\
&=& \sum_{i_{1} \neq \cdots \neq i_{r}} \sum_{q=1}^{r}
        x_{i_{1}} \cdots \,\, \widehat{x_{i_{q}}} \,\, \cdots x_{i_{r}}
        \sum_{j \neq i_{q}} \frac{1}{x_{j}-x_{i_{q}}}
\end{eqnarray*}
Write
\[ \sum_{j \neq i_{q}} \frac{1}{x_{j}-x_{i_{q}}}
    = \sum_{j \neq i_{1}, \cdots, i_{r}} \frac{1}{x_{j}-x_{i_{q}}}
        + \sum_{1 \leq p \leq r, \, p \neq q} \frac{1}{x_{i_{p}}-x_{i_{q}}}. \]
For fixed indices $i_{1}, \cdots, \widehat{i_{q}}, \cdots, i_{r}$,
\[ \sum \frac{1}{x_{j}-x_{i_{q}}} = 0 \]
 if the summation is running over all
possible values for $i_{q}$ and $j$ such that $i_{q} \neq j$ and
both of them are not equal to $i_{1}, \cdots, \widehat{i_{q}}, \cdots, i_{r}$.
Therefore
\begin{eqnarray}\label{eqn:dersrA}
\frac{r!}{m} \frac{d}{dt} y_{r}
&=& \sum_{i_{1} \neq \cdots \neq i_{r}}
        x_{i_{1}} \cdots  x_{i_{r}}
        \sum_{1 \leq p, q \leq r, \, p \neq q} \frac{1}{x_{i_{q}}(x_{i_{p}}-x_{i_{q}})}.
\end{eqnarray}
Write
\[ \sum_{1 \leq p, q \leq r, \, p \neq q} \frac{1}{x_{i_{q}}(x_{i_{p}}-x_{i_{q}})}
    = \sum_{p>q} \frac{1}{x_{i_{q}}(x_{i_{p}}-x_{i_{q}})}
        + \sum_{p<q} \frac{1}{x_{i_{q}}(x_{i_{p}}-x_{i_{q}})}. \]
Switch the indices $p$ and $q$ in the second term and adding the first term, we have
\begin{eqnarray*}
\sum_{1 \leq p, q \leq r, \, p \neq q} \frac{1}{x_{i_{q}}(x_{i_{p}}-x_{i_{q}})}
&=& \sum_{p>q} \frac{1}{x_{i_{p}}-x_{i_{q}}} \left(\frac{1}{x_{i_{q}}} - \frac{1}{x_{i_{p}}} \right) \\
&=&  \sum_{p>q} \frac{1}{x_{i_{p}} x_{i_{q}}}.
\end{eqnarray*}
Hence by equation \eqref{eqn:dersrA}
\begin{eqnarray*}
 \frac{r!}{m} \frac{d}{dt} y_{r}
&=&  \sum_{i_{1} \neq \cdots \neq i_{r}} \,\, \, \sum_{1 \leq p, q \leq r, \, p > q}
        x_{i_{1}} \cdots \,\, \widehat{x_{i_{q}}} \,\, \cdots \,\,
            \widehat{x_{i_{p}}} \,\,\cdots x_{i_{r}}  \\
&=& \frac{1}{2} r(r-1) (k-r+3)(k-r+2) (r-2)! y_{r-2}
\end{eqnarray*}
This proves the claim.

The explicit formula for $y_{r}(t)$ can be obtained from \eqref{ak1}
recursively, and  it is  a polynomial
in $t$ and initial conditions $x_{1}(0), \cdots, x_{k+1}(0)$.
For each $t$, we can obtain $x_{1}(t), \cdots, x_{k+1}(t)$ as the $k+1$ solutions
of the following polynomial equation in $z$:
\begin{equation} \label{eqn:rootssymAk}
 \sum_{r=0}^{k+1} \,\, (-1)^{k+1-r} \,\, y_{k+1-r}(t) \, z^{r} = 0,
\end{equation}
with the property
\[ x_{1}(t) < x_{2}(t) < \cdots < x_{k+1}(t). \]
\eeg

\beg \label{sec:Bk} {\bf The $B_{k}$ case}\par

\ss
Suppose that $k \geq 2$.
Let $\{ e_{1}, \cdots e_{k} \}$ be the standard orthonormal basis of $\mathbb{R}^{k}$
and $(x_{1}, \cdots , x_{k})$ the corresponding coordinate.
We identify $\mathbb{R}^{k}$ with a normal space of an isoparametric submanifold of type $B_{k}$.
The set  $e_{k}$ and
$\frac{1}{\sqrt{2}} \left(e_{i} - e_{i+1} \right)$ with $1 \leq i \leq k-1$ is a simple root system of $B_k$, and the set of positive roots are $e_{i}$ with $1 \leq i \leq k$ and
$\frac{1}{\sqrt{2}} \left(e_{i} \pm e_{j} \right)$ with $1 \leq i<j \leq k$.
The Coxeter group is generated by all permutations and sign changes of $\{e_{1}, \cdots, e_{k}\}$.
Since the multiplicities of curvature spheres are invariant under the action of the Coxeter group,
there are only two possible multiplicities.
Let $m_{1}$ be the multiplicities of curvature spheres corresponding to
$\frac{1}{\sqrt{2}} \left(e_{i} \pm e_{j} \right)$ , and $m_{2}$
the multiplicities of curvature spheres corresponding to $e_{i}$.  So the MCF \eqref{ab} is
$$-x'= m_1\left(\sum_{i<j} \frac{e_i+e_j}{x_i+x_j} + \frac{e_i-e_j}{x_i-x_j}\right)
+ m_2\sum_{i} \frac{e_i}{x_i}.$$
So
$$-x_i'= \frac{m_2}{x_i}+ m_1\left(\frac{1}{x_i-x_j}+ \frac{1}{x_i+x_j}\right).$$
Set $y_i= x_i^2$. Then we have
\begin{equation} \label{eqn:MCFB}
 \frac{1}{2} \frac{d}{dt} y_{i} =  - m_{2}
    - m_{1} \sum_{j \neq i} \frac{2 y_{i}}{y_{i} - y_{j}}.
\end{equation}
Let $s_0=1$,  $s_{i}$ the $i$-th elementary symmetric polynomial
of $y_1, \ldots, y_k$, and
$$\zeta_r(t)= s_r(x(t)).$$
We claim that
\beq\label{bk1}
\zeta_j'= -2(k-j+1)(m_2+m_1(k-j)) \zeta_{j-1},\quad 1\leq j\leq k.
\eeq
Note that when $j=1$ the right hand side is $-2k(m_2+m_1(k-1))$, which is equal to $-2n$.
To prove this claim, we compute as follows:
First
\begin{align*}
&\frac{1}{2} \frac{d}{dt} \left( y_{i_{1}} \cdots y_{i_{j}} \right)
= \sum_{l=1}^{j} \frac{y_{i_{1}} \cdots y_{i_{j}}}{y_{i_{l}}}
    \left( - m_{2} - 2 m_{1} \sum_{p \neq i_{l}} \frac{y_{i_{l}}}{y_{i_{l}} - y_{p}} \right) \\
&\,\, =  y_{i_{1}} \cdots y_{i_{j}}
    \left( - m_{2} \sum_{l=1}^{j} \frac{1}{y_{i_{l}}}
    - 2 m_{1} \sum_{l=1}^{j} \sum_{p \neq i_{l}} \frac{1}{y_{i_{l}}-y_{p}}
    \right).
\end{align*}
Since
\begin{eqnarray*}
\sum_{l=1}^{j} \sum_{p \neq i_{l}} \frac{1}{y_{i_{l}}-y_{p}}
&=& \sum_{l=1}^{j} \sum_{1 \leq q \leq j, \,\, q \neq l} \frac{1}{y_{i_{l}}-y_{i_{q}}}
    + \sum_{l=1}^{j} \sum_{p \neq i_{1}, \cdots, i_{j}} \frac{1}{y_{i_{l}}-y_{p}},
\end{eqnarray*}
and the first term on the right hand side is 0, we have
\begin{align}
&\frac{d}{dt} s_{j}
= \frac{1}{j!} \frac{d}{dt} \left(\sum_{i_{1} \neq \cdots \neq i_{j}} y_{i_{1}} \cdots y_{i_{j}}
\right)  \nonumber \\
&\,\,= - \frac{2 m_{2}}{j!}  \sum_{l=1}^{j} \sum_{i_{1} \neq \cdots \neq i_{j}}
    \frac{y_{i_{1}} \cdots y_{i_{j}}}{y_{i_{l}}}
    - \frac{4 m_{1}}{j!}\sum_{l=1}^{j}  \sum_{i_{1} \neq \cdots \neq i_{j} \neq p}
    \frac{y_{i_{1}} \cdots y_{i_{j}}}{y_{i_{l}}-y_{p}}. \label{eqn:dersigma}
\end{align}
The first term on the right hand side of this equation is
$ -2 m_{2} (k-j+1) \zeta_{j-1}$.
 To compute the second term, we notice that
\[\sum_{i_{1} \neq \cdots \neq i_{j} \neq p}
    \frac{y_{i_{1}} \cdots y_{i_{j}}}{y_{i_{l}}-y_{p}}
    \]
can be written as
\[
 \sum_{ i_{1} \neq \cdots \neq i_{j}, i_{l} > p }
       \frac{( y_{i_{1}} \cdots \widehat{y_{i_{l}}} \cdots y_{i_{j}})y_{i_{l}}}{y_{i_{l}}-y_{p}}
    + \sum_{ i_{1} \neq \cdots \neq i_{j},  i_{l} < p}
    \frac{ (y_{i_{1}} \cdots \widehat{ y_{i_{l}}} \cdots y_{i_{j}})y_{i_{l}}}{y_{i_{l}}-y_{p}}.
\]
Switching the indices $i_{l}$ and $p$ in the second term and adding to the first term, we have
\begin{eqnarray*}
\sum_{i_{1} \neq \cdots \neq i_{j} \neq p}
    \frac{y_{i_{1}} \cdots y_{i_{j}}}{y_{i_{l}}-y_{p}}
&=& \sum_{\begin{array}{c} i_{1} \neq \cdots \neq i_{j} \neq p \\ i_{l} > p \end{array}}
    y_{i_{1}} \cdots \widehat{y_{i_{l}}} \cdots y_{i_{j}}
\end{eqnarray*}
Therefore the second term on the right hand side of equation \eqref{eqn:dersigma} is
\[ -2 m_{1}(k-j+1)(k-j) s_{j-1}.\]
This proves the claim.

Note that explicit formula for $\zeta_{i}(t)$ can be obtained from \eqref{bk1} recursively, and
it is always a
degree $i$ polynomial
in $t$ and initial conditions $y_{1}(0), \cdots, y_{k}(0)$.
 Let $y_{1}(t), \cdots y_{k}(t)$
be the $k$ roots of
\begin{equation} \label{eqn:rootssymBk}
 \sum_{r=0}^{k} \,\, (-1)^{k-r} \,\, s_{k-r}(t) \, z^{r} = 0,
\end{equation}
with the property $y_{1}(t) > y_{2}(t) > \cdots > y_{k}(t) > 0$.
Then $x_{i}(t) = - \sqrt{y_{i}(t)}$ for $i=1, \cdots k$.
\eeg

\beg \label{sec:Dk} {The $D_{k}$ case}\par

\ss
Suppose that $k \geq 4$.
Let $\{ e_{1}, \cdots e_{k} \}$ be the standard orthonormal basis of $\mathbb{R}^{k}$
and $(x_{1}, \cdots , x_{k})$ the corresponding coordinate.
We identify $\mathbb{R}^{k}$ with a normal space of an isoparametric submanifold of type $D_{k}$.
The set of simple roots are $\frac{1}{\sqrt{2}} \left(e_{k-1} + e_{k} \right)$ and
$\frac{1}{\sqrt{2}} \left(e_{i} - e_{i+1} \right)$ with $1 \leq i \leq k-1$, and the set
of positive roots is
$\{\frac{1}{\sqrt{2}} \left(e_{i} \pm e_{j} \right)\n 1 \leq i<j \leq k\}$.
The Coxeter group is generated by all permutations and even number of sign changes of
$\{e_{1}, \cdots, e_{k}\}$.
The open positive Weyl chamber $C$ is
\[ C := \{ x \in \mathbb{R}^{k} \mid x_{1} < x_{2} < \cdots < x_{k} \,\,\,
    {\rm and} \,\,\,  x_{k-1}+x_{k} < 0 \}. \]
All multiplicity for curvature spheres are equal and will be denoted by $m$.
The mean curvature flow equation \eqref{ab} becomes
\beq\label{dk1}
-x'= m\, \sum_{i<j} \frac{e_i-e_j}{x_i-x_j} + \frac{e_i+e_j}{x_i+x_j}.
\eeq
\[ \frac{1}{m} \frac{d}{dt} x_{i}
    = \sum_{q \neq i} \frac{2 x_{i}}{x_{q}^{2}-x_{i}^{2}}. \]
Multiply both sides by $\frac{x_{i}}{2}$ to obtain
\begin{equation} \label{eqn:MCFD}
 \frac{1}{4m} \frac{d}{dt} y_{i} =  \sum_{q \neq i} \frac{ y_{i}}{y_{q} - y_{i}}
\end{equation}
where
$y_{i} :=  x_{i}^{2}$ for all $i=1, \cdots, k$.

Let $s_{i}$ be the $i$-th elementary symmetric polynomial
of $y_{1}, \cdots, y_{k}$ and we set $s_{0}=1$.
Then $s_{1}, \cdots, s_{k-1}$ and $\sqrt{s_{k}}$
generate the algebra of polynomials invariant under the action of the Coxeter
group.
Note that equation \eqref{eqn:MCFD} is a special case of the equation
\eqref{eqn:MCFB} with $m_{2}=0$ and $m_{1} = m$.
Set  $\zeta_j(t)= s_j(x(t))$.
The proof of \eqref{bk1} also works here, and we obtain
$$\zeta_{r}' = -2 m (k-r+1)(k-r)\zeta_{r-1}, \qquad 1\leq r\leq k.$$
\eeg

\ms
\beg {\bf The rank $2$ cases}\par

\ss
We use a different set of generators for the ring of $W$-invariant polynomials to compute explicit
solutions of the MCF \eqref{aa} for the rank $2$ cases.  The Weyl group is the dihedral group with
$2g$ elements. We identify the normal space of a rank $2$ isoparametric submanifold with $\R^2=\C$,
and use $e^{\frac{ik\pi}{g}}$ with $0\leq k< g$ as positive roots.  Then
$$P_1(x)= x_1^2+x_2^2, \quad P_2(x)= {\rm Re}((x_1+ix_2)^g)$$
form a set of generator of the ring of invariant polynomials
($g=3$ for $A_2$, $g=4$ for $B_2$, and $g=6$ for $G_2$). It is known (cf. [PTb]) that
\ben
\item if  $g=3, 6$, then all multiplicities are equal and are either $1$ or $2$; and if $g=4$, then
there are two positive integers $m_1\leq m_2$ such that the multiplicity corresponding
to $\R \vn_j= \R e^{\frac{j \pi i}{4}}$ is $m_1$ for $j$ even and is $m_2$ for $j$ odd.
\item Let $F_i= \sum_{i=1}^g \frac{m_i\li \K P_i, \vn_i\ri}{\li x, \vn_i\ri}$ for $i=1, 2$ be as in Theorem \ref{ba}. Then
$$F_1(x)=2n,\quad F_2(x)=\bca 0 &\, {\rm if\,\, } g=3 \,\, {\rm or \,\, 6 \/},
\\ 8(m_2-m_1)(x_1^2+x_2^2), & \,{\rm if\,\,} g=4,  \eca$$
\een
If $x(t)$ is a solution of the MCF, then it follows from Theorem \ref{ba} that $y(t) =(P_1(x(t)), P_2(x(t))$ satisfies
\beq\label{ah}
\bca y_1'= -2n,\\
y_2'=\bca 0, & {\rm if\,} g=3, 6,\\
-8(m_2-m_1) y_1, & {\rm if\, } g=4.
\eca\eca
\eeq
By Lemma \ref{ac}, it suffices to consider initial data $x_0=e^{i\o_0}$ for the MCF.
The corresponding solution for \eqref{ah} with initial data $y(0)=(1, \cos g \o_0)$ is
\begin{align*}
& y_1(t)= 1-2nt, \quad y_2(t)= \cos(g\o_0), \quad {\rm if\,\,} g= 3, 6,\\
& y_1(t)=1-2nt, \quad y_2(t)= \cos 4\o_0 - 8(m_2-m_1)(t-nt^2), \quad {\rm if\,\,} g=4.
\end{align*}
  Set $x(t)= r(t)e^{i\o(t)}$. Since $y_1(t)= r^2(t)$ and $y_2(t)= r^g(t) \cos g\o(t)$,
  solution of the MCF for rank $2$ case with initial data $e^{i\o_0}$  is
$r(t)= (1-2nt)^{\frac{1}{2}}$ and
$$ \o(t)=\bca \frac{1}{g}\,\cos^{-1}\left(\frac{\cos g\o_0}{(1-2nt)^{\frac{g}{2}}}\right), & g= 3, 6,\\
\frac{1}{4}\, \cos^{-1} \left(\frac{\cos 4\o_0 - 8(m_2-m_1) (t-nt^2)}{(1-2nt)^2}\right), & g=4.
\eca$$

We claim that the isoparametric submanifold through $x_0= e^{i\o_g}$ is minimal in sphere, where
\begin{equation}\label{bb}
\o_g= \bca \frac{\pi}{2g}, &  g=3, 6,\\ \frac{1}{4}\cos^{-1}(\frac{m_2-m_1}{m_2+m_1}), & g=4.\eca
\end{equation}
To see this,   by Lemma \ref{miniso}, the the polar angle of the MCF flow starting at a minimal
submanifold in sphere must be some constant $\o_0$.  So:
\ben
\item[(i)] For $g=3, 6$, we have $y_2(t)= r(t)^g \cos g\o_0= \cos g\o_0$.
This implies that $\cos g\o_0= 0$, hence $\o_0= \frac{\pi}{2g}$.
\item[(ii)] For $g=4$,  $y_2(t)= y_1^2(t) \cos 4\o_0$ implies that
$$y_2'= 2y_1y_1'\cos 4\o_0= -4n y_1\cos 4\o_0= -8(m_2-m_1)y_1.$$
Hence $\cos 4\o_0=\frac{2(m_2-m_1)}{n}=\frac{m_2-m_1}{m_1+m_2}$.
This proves the claim.
\een
\eeg

Next we compute the maximal time $T$ for the MCF for the above examples.
The flow blows up  when $\cos g\o(t)= \pm 1$.
If $g=3$ or $6$, then
$$y_2(t)= r(t)^g \cos g\o(t) = (1-2nt)^{g/2} \cos g\o(t)= \cos g\o_0.$$
Hence  $T= \frac{1}{2\pi} (1-|\cos g\o_0)|^{2/g})$.
For $g=4$, we have
$$y_2(t)= r(t)^4 \cos 4\o(t) = (1-2nt)^2 \cos 4\o(t) = \cos 4\o_0 - 8(m_2-m_1)(t-nt^2).$$
Then $\cos g\o(T)=1$ if $\o_0\in (0, \o_4)$, and $\cos g\o(T)=-1$ if $\o_0\in (\o_4, \frac{\pi}{4})$.
This solves $T$ and we get:

\bprop Let $\o_g$ be the constant defined by \eqref{bb}.  Then
\ben
\item The MCF through $e^{i\o_g}$  homothetically converges to $0$ at time $T=\frac{1}{2n}$.
\item For $g= 3, 6$, the maximal time for the MCF with initial data $e^{i\o_0}$ is
$T= \frac{1}{2n}(1-|\cos g\o_0|^{\frac{2}{g}})$. For $g=4$, the maximal time is
\[ T =\left\{ \begin{array}{ll}
    \frac{1}{2n} \left(1-\sqrt{\frac{m_{1}+m_{2}}{2 m_{1}}
                (\cos 4\o_0- \frac{m_2-m_1}{m_2+m_1})}\right),
        & {\rm if} \,\,\, \theta_{0} \in (0, \theta_{4}), \\ \\
    \frac{1}{2n} \left(1-\sqrt{\frac{m_{1}+m_{2}}{2 m_{2}}
                (-\cos 4\o_0+ \frac{m_2-m_1}{m_2+m_1})}\right),
        & {\rm if} \,\,\, \theta_{0} \in (\theta_{4}, \frac{\pi}{4}).
        \end{array} \right. \]
\item
$$\lim_{t\to T^-} \o(t)= \bca 0, &{\rm if \,\, } \o_0\in (0, \o_g),\\
        \frac{\pi}{g}, & {\rm if\,\,} \o_0\in (\o_g, \frac{\pi}{g}).\eca$$
\een
\eprop

\bs

\section{Mean curvature flows of focal submanifolds}
\label{sec:focal}

We consider the mean curvature flows
of focal submanifolds of isoparametric submanifolds. They behave very similarly to the
mean curvature flows of isoparametric submanifolds. With slight modification,
most results in section \ref{sec:Basic}
also hold for mean curvature flows of focal submanifolds.

Let $M_{0}$ be an isoparametric submanifold and $p \in M_{0}$ a fixed point.
As before, we assume that $M_{0}$ is contained in a sphere centered at the origin.
Let $C \subset p+ \nu_{p}M_{0}$ be the Weyl chamber containing $p$, and
$\Delta^{+}$ the set of positive roots of the Weyl group $W$ associated to $M_0$.
 Then
\ben
\item $\bar  C$ is a stratified space,
\item the isotropy subgroup of any two points in the same stratum $\sigma$ are the same,
and will be denoted by $W_\sigma$,
\item for $x\in \p C$, let
$$\D^+(x)=\{\a\in \D^+\n \li x, \a\ri >0\},$$ then $\D^+(x_1)=\D^+(x_2)$ if and only if
$x_1, x_2$ lie in the same stratum $\sigma$, and will be denoted by $\D^+(\sigma)$,
\item $\sigma$ is the Weyl chamber of $W_x$ for $x\in \sigma$ and $\sigma$ is an
open simplicial cone in the following  linear subspace
$$V(\sigma)=\{x\in p+\nu_pM_0\n \li x, \a\ri
=0, \,\, {\rm for\, all\,} \, \a\in \D^+\setminus \D^+(\sigma)\},$$
\een

Let $\sigma$ be a stratum in $\p C$, $x_0\in \sigma$, and $M$ the focal submanifold of
$M_0$ through $x_0$.
By \cite[Theorem 4.1]{Terng},
the mean curvature vector field of $M$ at $x_0$ is given by
\begin{equation} \label{eqn:focalMCV}
H(x_0) = - \sum_{\alpha \in \Delta^{+}(\sigma)} \frac{m_{\alpha}}{\li x_0, \alpha\ri} \,\, \alpha
\end{equation}
where $m_{\alpha}$ are multiplicities of curvature spheres of $M_{0}$.
Moreover
\begin{equation} \label{eqn:focalxH}
 \li x_0, \alpha\ri \,\, =\,\, \li H(x_0), \alpha\ri \,\,= \,\,0
\end{equation}
for all $\alpha \in \Delta^{+} \setminus \Delta^{+}(\sigma)$.
The mean curvature flow equation of $M$ is the following ODE on $\sigma$:
\begin{equation}  \label{eqn:focalMCF}
\frac{dx}{dt}  = - \sum_{\alpha \in \Delta^{+}(\sigma)}
\frac{m_{\alpha}}{\li x, \alpha\ri} \,\, \alpha.
\end{equation}

The analogue of Theorem \ref{thm:radialMCV} also holds for this case. In particular,
if $x(t)$ satisfies the flow equation \eqref{eqn:focalMCF} then
\begin{equation}
 \|x(t)\|^{2} = \|x(0)\|^{2} - 2nt
\end{equation}
where $n = \sum_{\alpha \in \Delta^{+}(\sigma)} m_{\alpha}$ is the dimension of $M$.
Therefore we have
\begin{thm}
The maximal interval for the solution of the mean curvature flow equation for any focal submanifold
is finite.
\end{thm}

Suppose $x(t)$ and $y(t)$ satisfy equation \eqref{eqn:focalMCF} on $\sigma$ and $x(0)\neq y(0)$. Use the same computation for \eqref{bp} to get
\begin{align}\label{bd}
\frac{1}{2} \frac{d}{dt} \|x(t) - y(t) \|^{2}
&= \li x(t) - y(t),  \,\, x^{\prime}(t) - y^{\prime}(t)\ri   \nonumber \\
&= \sum_{\a\in \D_+(\sigma)}
    m_{\a}\frac{\li x(t)-y(t), \a \ri^2}{\li x(t), \a \ri \li y(t), \a \ri}\, > 0.
    \end{align}
Then proofs given in section \ref{sec:Basic} works, so we have

\bthm \label{bc}
Let $M^n\subset S^{n+k-1}$ be a compact isoparametric submanifold in $\R^{n+k}$, $W$ its Weyl group,
$C$ the Weyl chamber in $x_0+\nu(M)_{x_0}$ containing $x_0\in M$, and $M_y$ the submanifold parallel
to $M$ through $y$.  If $\sigma \subset \overline{C}$ is a stratum, then
\ben
\item there is a unique $y_\sigma\in \sigma$ such that the focal submanifold $M_{y_\sigma}$ is minimal
in $S^{n+k-1}$, and the MCF in $\R^{n+k}$ with initial data $M_{y_{\sigma}}$ homothetically shrinks
to a point,
\item if $y_0\in \sigma \bigcap S^{n+k-1}-\{y_\sigma\}$, then the MCF in $\R^{n+k}$ with $M_{y_0}$
as initial data blows up in finite time $T < \frac{1}{2n}$, $x(t)\in \sigma$ for all $t\in [0, T)$,
and
$\lim_{t\to T^-} x(t)\,\in \p \sigma$, in particular, the limit is a focal submanifold with
lower dimension,
\item if $y_1, y_2\in \sigma \bigcap S^{n+k-1}$ are distinct, then
the MCF in $\R^{n+k}$ with initial data $M_{y_1}$ and $M_{y_2}$  converge to
distinct focal submanifolds of lower dimensions.
\een
\ethm

\bs

\section{Mean curvature flows for isoparametric submanifolds in spheres}
\label{sec:sphere}

If $M^{n}\subset S^{n+k-1}$ is an isoparametric submanifold in   $\R^{n+k}$, then
$M$ is also isoparametric
in $\mathbb{R}^{n+k}$. So basic structure theory for isoparametric submanifolds in Euclidean
spaces also applies to $M$.
For $x \in M$, let $H(x)$ and $H_{E}(x)$ be the mean curvature vector fields of $M$ at $x$
as a submanifold of $S^{n+k-1}$ and $\mathbb{R}^{n+k}$ respectively.
Then $H(x)$ is the orthogonal projection of $H_{E}(x)$ to $T_{x}S^{n+k-1}$. More precisely
\[ H(x) = H_{E}(x) + n x \]
for all $x \in M$. In particular, $H$ is again a parallel normal vector field along $M$.
The mean curvature flow of $M$ as a submanifold of $S^{n+k-1}$ behaves similarly
to its flow as a submanifold of $\mathbb{R}^{n+k}$. With slight modifications, most results
for mean curvature flows for isoparametric submanifolds in the Euclidean spaces also hold
for isoparametric submanifolds in spheres. We only need to explain how to deal with
the arguments in the Euclidean case which can not be applied directly to the spherical case.

Fix $x_{0} \in M$ and let $V = x_0+ \nu_{x_{0}}M$ be
the normal space of $M$ as a submanifold of $\mathbb{R}^{N}$ at the point $x_{0}$, $W$ its
Coxeter group, and $C\subset V$ the Weyl chamber containing $x_0$.
The mean curvature flow of $M$ in $S^{n+k-1}$ is uniquely determined by the flow of $x_{0}$ in
$S := C \bigcap S^{k-1}$:
\begin{equation} \label{be}
x'(t)= -\sum_{\a\in \D_+} \frac{m_\a \a}{\li x(t),\a\ri} \, + n x(t).
\end{equation}
The set $S$ is a geodesic $(k-1)$-simplex on $S^{k-1}$.
Let $x(t) \in S$ be a solution to equation \eqref{be} with initial condition $x_{0}$.
Then
\[ y(t) = \sqrt{1-2nt} \,\,\,\, x\left(- \frac{1}{2n} \log(1-2nt)\right) \]
satisfies the Euclidean mean curvature flow equation \eqref{ab} with initial condition
$y(0) = x_{0}$.
Let $[0, T_{x})$ and $[0, T_{y})$ be the maximal intervals for the domains of $x(t)$ and $y(t)$
respectively. Then
\[ T_{x} = - \frac{1}{2n} \log(1-2n T_{y}) \]
and
\[ \lim_{t \rightarrow T_{y}^{-}} \,\, y(t)
    \,\,=\,\, \sqrt{1-2n T_{y}}\,\,  \lim_{t \rightarrow T_{x}^{-}} \,\, x(t). \]
Note that by Theorems \ref{thm:radialMCV} and Corollary \ref{cor:miniso},
$T_{y} \leq \frac{1}{2n}$ and the equality holds if and only if
the isoparametric submanifold $M_{0}$ passing $x_{0}$ is minimal in the sphere $S^{n+k-1}$.
So by Theorem \ref{thm:finiteconv}, if $M_{0}$ is not minimal in the sphere, then
$x(t)$ converges to a focal submanifold at a finite time $T_{x}$.
This proves Theorem \ref{thm:MCFisopS}.

If $x_{1}(t) \in S$ and $x_{2}(t) \in S$ satisfy the spherical mean curvature flow equation
\eqref{be}, then
\begin{eqnarray*}
\frac{1}{2} \frac{d}{dt} \|x_{1}(t) - x_{2}(t)\|^{2}
&=& \li x_{1}-x_{2}, \left(H_{E}(x_{1}) + n x_{1} \right)
- \left(H_{E}(x_{2}) + n x_{2} \right)\ri \\
&=& n \|x_{1}-x_{2}\|^{2} + \li x_{1}-x_{2}, H_{E}(x_{1})  - H_{E}(x_{2})\ri.
\end{eqnarray*}
By \eqref{bp}, $\li x_{1}-x_{2}, H_{E}(x_{1}) - H_{E}(x_{2})\ri \,\,\geq \,\, 0$.
Therefore
\begin{equation} \label{eqn:departS}
 \frac{d}{dt} \|x_{1}(t) - x_{2}(t)\|^{2} \,\,\geq \,\, 2n \|x_{1}(t)-x_{2}(t)\|^{2}.
\end{equation}
We use \eqref{eqn:departS} to give an estimate of the maximal interval $[0, T)$
for the spherical mean curvature flow $x(t)$.
Let $p_{0}$ be the unique point in $S$ such that the isoparametric submanifold
passing $p_{0}$ is minimal in the sphere $S^{n+k-1}$.
Set $x_{1}(t) = x(t)$ and $x_{2}(t) = p_{0}$ in equation \eqref{eqn:departS}. Since $x_2(t)$ exists for all $t>0$, we obtain
\[ \|x(t) - p_{0} \| \geq e^{nt} \|x(0) - p_{0} \| \]
 for all $t$ as long as $x(t)\in S$.
Let $D$ be the diameter of $S$, then $D \leq 2$ and
\[ T \leq \frac{1}{n} \log \frac{D}{\|x(0) - p_{0} \|}.\]

Now we discuss the behavior of invariant polynomials under the spherical mean curvature flow.
Let $x(t) \in S$ be the mean curvature flow of $x_{0}$.
For any function $f$ on $V$, let $f(t) = f(x(t))$. Then
\[ f^{\prime}(t) \,\,=\,\, \li \nabla f(x(t)), H(x(t))\ri \,\,=\,\,
    \li \nabla f(x(t)), H_{E}(x(t))\ri + n \li \nabla f(x(t)), x(t)\ri. \]
If $f$ is a homogenous polynomial of degree $k$ which is invariant under the action
of the Coxeter group $W$, then as in the proof of Theorem \ref{thm:finiteconv},
\begin{equation} \label{eqn:invpolES}
 f^{\prime}(t) \,\,=\,\, - F(x(t)) + nk f(t)
\end{equation}
where $F$ is defined by equation \eqref{eqn:mcfpoly}
and it is an invariant polynomial of degree $k-2$.
If we have computed $F(t) := F(x(t))$, then we can solve $f(t)$ from equation \eqref{eqn:invpolES}
and obtain
\begin{equation} \label{eqn:intinvpolS}
 f(t) = - e^{knt} \int e^{-knt} F(t) \, dt.
\end{equation}
Note that there is no homogeneous invariant polynomial of degree 1. By induction on the degree,
we obtain the following
\begin{thm}
If $x(t)$ satisfies the spherical mean curvature flow equation \eqref{be} and $f$ is a
$W$-invariant polynomial, then $f(t) = f(x(t))= c_1 e^{knt} +h(t)$ for some constant $c_1$  and polynomial $h$.
\end{thm}

In particular $f(t)$ is well defined for all $t \in \mathbb{R}$.
In section \ref{sec:solution}, we have given  explicit formulas for
$F_i$ for invariant homogeneous polynomials $P_i$ for isoparametric submanifolds.
We can use these formula and \eqref{eqn:intinvpolS} to construct  explicit
solutions to the spherical mean curvature flow equation for isoparametric submanifolds in spheres.

\beg {\bf Phase portrait for rank $2$ cases}\par

Let $M^n\subset S^{n+1}\subset \R^{n+2}$ be an isoparametric hypersurface with $g$ distinct principal curvatures.  Then the Weyl group associated to $M$ as a rank $2$ isoparametric submanifold in $\R^{n+2}$ is the dihedral group of $2g$ elements. Let $C$ denote the Weyl chamber containing $x_0\in M$, and $D$ the intersection of $C$ and the normal circle at $x_0$ in $S^{n+1}$.  Let $p_1, p_2$ denote the end points of $D$.  The  arc $D=\widehat{p_1p_2}$ has length $\pi/g$. For $y\in \bar C$, let $M_y$ denote the submanifold through $y$ that is parallel to $M$ (a leaf of the isoparametric foliation).  There exists a unique $p_0\in D$ such that $M_{p_0}$ is minimal in $S^{n+1}$.
\ben
\item The spherical MCF \eqref{be} has three orbits: a stationary point $p_0$, the orbit $\widehat{p_0p_1}$ with one end tends to $p_0$ and the other end tends to $p_1$, and the orbit $\widehat{p_0p_2}$ with one end tends to $p_0$ and the other end tends to $p_2$.
\item The MCF \eqref{aa} in $\R^{n+2}$ starting at $M_y$ degenerates homothetically to one point (the origin) if $y=p_0$, to $M_{rp_2}$ for some $0<r<1$ if $y\in \widehat{p_0p_2}$, and to $M_{rp_1}$ for some $1<r<1$ if $y\in \widehat{p_1p_0}$.
\een
\eeg

\beg {\bf Phase portrait for the $A_3$ cases}\par

\ss
Let $M^n\subset S^{n+2}$ be an isoparametric submanifold with Weyl group $A_3$ and uniform
multiplicity $m$, and $x_0\in M$. Let $C$ denote the Weyl chamber containing $x_0$, and $D$
the intersection of $C$ and the normal sphere at $x_0$.  Then $D$ is a geodesic triangle
with  vertices $p_1, p_2, p_3$ and interior angles $\frac{\pi}{3}, \frac{\pi}{3}, \frac{\pi}{2}$.
The phase spaces of spherical MCF \eqref{be} and Euclidean MCF \eqref{aa} are $D$ and $C$ respectively.   We describe
the phase portraits:

(1) There exists a unique $p_0$ in the interior of $D$ such that $p_0$ is the fixed point
of the ODE \eqref{be}.  This implies that the spherical MCF starting at the parallel
submanifold $M_{p_0}$ is stationary and $M_{p_0}$ is minimal in $S^{n+k-1}$.

(2) For each $1\leq i\leq 3$, there is a unique flow $\ell_i$ that
at one end approaches  $p_0$ and at the other end approaches
$p_i$. This implies that for $y\in \ell_i$, the spherical MCF
starting at $M_y$ collapses in finite time to the focal
submanifold $M_{p_i}$ by collapsing fibers of the fibration
$M_y\to M_{p_i}$. The fibers are isoparametric submanifolds with
Weyl group $A_2$ for $i=1, 2$ and $A_1\times A_1$ for $i=3$.  In
particular, when $m=2$, $M_y$ is diffeomorphic to the manifold of
flags in $\C^4$ and collapsing is along complex flag manifolds of
$\C^3$ in $M_y$ for $i=1, 2$ and $S^2\times S^2$ in $M_y$ for
$i=3$.

(3) For distinct $i, j, k$, let $D_k$ denote the triangle  with vertices $p_i, p_j, p_0$ and
edges $\ell_i, \ell_j$, and geodesic segment $\widehat{p_ip_j}$ in the sphere. The flow for \eqref{be} starting
at a point in the interior $D_k^0$ of $D_k$ exists for finite time and converges to a point on
the interior of $\widehat{p_ip_j}$.  This implies that for $y\in D_k^0$, the spherical MCF
starting at $M_y$ converges in finite time to a focal submanifold $M_q$ with
$q\in\widehat{p_ip_j}\setminus \{p_i, p_j\}$ by collapsing one family of curvature spheres.

(4) The flow of \eqref{aa} on $C$ starting at $p_0$ is the straight line joining the
origin to $p_0$, the flow starting from a point in $D_k^0$ converges to a point on the
wall containing $p_i, p_j$ ($i, j, k$ distinct), and the flow starting at a point on $\ell_i$
converges to a point on the line segment $\overline{Op_i}$.   This implies that the Euclidean
MCF with initial data $M_y$
\ben
\item[(i)] shrinks homothetically to the origin if $y=p_0$,
\item[(ii)] converges to a focal submanifold $M_q$ for some $q$ lies in the open cone spanned
by $\widehat{p_ip_j}$ and the collapsing is along a curvature $m$-sphere if $y\in D_k^0$,
\item[(iii)] converges to a focal submanifold $M_q$ with $q\in \overline{Op_i}$ for $y\in \ell_i$,
moreover, the collapsing is along fibers of the fibration $M_y\to M_q$ and the fibers are
isoparametric submanifolds with Weyl group $A_2, A_2$, and $A_1\times A_1$ respectively for
$i=1, 2, 3$.
\een
\eeg

\bs
\section{Mean curvature flow for polar action orbits}
\label{sec:polar}

Let $G$ act on a Riemannian manifold $N$ isometrically, and $G\cdot p$ be a principal orbit
through $p$.  If $v\in \nu(G\cdot p)_p$, then $\hat v(g\cdot p)= dg_p(v)$ is a globally defined
normal vector field on $G\cdot p$ and is called a {\it $G$-equivariant normal field}.
An isometric action of a compact Lie group $G$ on a Riemannian manifold $N$ is called
{\it polar\/} if there is a totally geodesic submanifold $\Sigma$ that meets all $G$-orbits
and meets orthogonally.  Such $\Sigma$ is called a {\it section\/}.  We list some properties
of polar actions (cf. \cite{PT}):
\ben
\item The action $G$ is polar if and only if every $G$-equivariant normal field is parallel
with respect to the induced normal connection on $G\cdot p$ as a submanifold of $N$.
\item Let $N(\Sigma)=\{g\in G\n g\cdot \Sigma=\Sigma\}$ and
$Z(\Sigma)=\{ g\in G\n g\cdot x= x \,\, \forall\,\, x\in \Sigma\}$ denote the normalizer and
centralizer of $\Sigma$ respectively. Then the quotient group $W(\Sigma)=N(\Sigma)/Z(\Sigma)$ is a finite
group acting on $\Sigma$, and is called {\it the generalized Weyl group associated to the
polar action\/}.
\item The orbit space $\Sigma/W$ is isomorphic to $N/G$ and the ring of smooth $G$-invariant
functions on $N$ is isomorphic to the ring of $W$-invariant functions on $\Sigma$ under the
restriction map.
\item If $p_0\in \Sigma$ is a singular point, i.e., $G\cdot p_0$ is a singular orbit, then
the slice representation of $G_{p_0}$ on the normal space of the orbit at $p_0$ is a polar
representation.
\een

\bthm
Suppose the isometric action $G$ on $N$ is polar, and $\Sigma$ is a section.  Then
\ben
 \item[(i)] if $x\in \Sigma$, then the mean curvature vector $\xi(x)$ of $G\cdot x$ at $x$
 is tangent to $\Sigma$ at $x$,
 \item[(ii)] if $x'(t)= \xi(x(t))$ with $x(t)$ regular (i.e., $G\cdot x(t)$ is a principal orbit),
 then $G\cdot x(t)$ satisfies the MCF in $N$, in other words, the MCF in $N$ with a principal
 $G$-orbit as initial data flows among principal $G$-orbits.
 \een
\ethm

For general polar action, $W$ need not be a Coxeter group and the orbit space of the
$W$-action on the section can be complicated. In fact,
given any finite group $W$ and any compact $W$-manifold, there exist a
Riemannian manifold $N$, a compact group $G$, and an isometric  polar $G$-action on $N$  such that the induced action on the section is the given $W$-action (cf. \cite{PT}).
 Hence the behavior of the MCF for general polar
actions is not as clear as in the sphere and Euclidean case.

A polar action on a symmetric space is {\it hyperpolar\/} if the sections are flat. In this case
the fundamental domain of the $W$-action on a section is a geodesic simplex.  A submanifold in a symmetric space is called {\it equifocal\/} if its normal bundle is flat, exponential of each normal space is a flat, and the focal radii along a parallel normal field are constant.  It was proved in \cite {TT} that principal orbits of a hyperpolar action on symmetric space are equifocal, and  parallel foliation of an equifocal submanifold is an orbit-like foliation. Moreover, generators of the ring of smooth functions that are constant along parallel leaves were constructed in \cite{HLO}.    Hence we believe that methods developed in this paper can be used to solve the MCF starting with an equifocal submanifold  in symmetric spaces.

%%%%%%%%%%%%%%%%%%%%%%%%%%%%%%%%%%%%%%%%%%%%%%%%%%%


\begin{thebibliography}{399}
\bibitem[FKM]{FKM}D. Ferus, H. Karcher, and H.F. M\"{u}nzner,
{\it Cliffordalgebren und neue isoparametrische hyperflachen}, Math. Z. 177 (1981), 479-502.
\bibitem[GH]{GH} M.E. Gage and R.S. Hamilton, {\it The heat equation shrinking convex plane curves}, J. Differential Geometry 23 (1986), 69-96.
\bibitem[GB]{GB} L.C. Grove and C.T. Benson,
{\it Finite reflection groups}, second edition, GTM 99, Springer-Verlag, 1985.
\bibitem[HLO]{HLO} E. Heintze, X. Liu, C. Olmos,
{\it Isoparametric submanifolds and a Chevalley-type restriction theorem}, Integrable systems, geometry, and topology, AMS/IP Stud. Adv. Math., 36 (2006), 151--190.
\bibitem[HOT]{HOT} E. Heintze, C. Olmos, and G. Thorbergsson,
{\it Submanifolds with constant principal curvatures and normal holonomy group},
Int. J. Math. 2 (1991), 167-175.
\bibitem[Hu]{Hu} G. Huisken, {\it Contracting convex hypersurfaces in Riemannian manifolds by their mean curvature}, Invent. Math. 84 (1986), 463-480.
\bibitem[PT]{PT}
R.S. Palais, C.-L. Terng, {\it A general theory of canonical forms},  Transaction, Amer.
Math.   Soc. 300 (1987) 771-789.
\bibitem[PTb]{PTb}
R.S. Palais, C.-L. Terng, {\it Critical Point Theory and Submanifold Geometry},
Lecture Notes in Math., vol. \textbf{1353} (1988),  Springer-Verlag
\bibitem[T]{Terng} C.-L. Terng,
{\it Isoparametric submanifolds and their Coxeter groups},
J. Differential Geometry, 21 (1985), 79-107.
\bibitem[TT]{TT}  C.-L. Terng, G. Thorbergsson,
{\it Submanifold Geomtry in Symmetric spaces},
J. Differential Geometry, 42 (1995) 665-718.
\bibitem[Th]{Th} G. Thorbergsson,
{\it Isoparametric foliations and their buildings}, Ann. of Math., 133 91991), 429-446.
\bibitem[V]{Var} V.S. Varadarajan,
{\it Lie groups, Lie algebras, and their representations}, GTM 102, Springer-Verlag, 1984.
\bibitem[W]{W} M.-T. Wang,
{\it Mean curvature flow in higher codimension}, to appear in the Proceedings of International Congress of Chinese Mathematicians 2001,
math.DG/0204054
\bibitem[Z]{Zhu} X.-P. Zhu,
{\it Lectures on mean curvature flows}, Studies in Advanced Math., AMS/IP, 2002.
\end{thebibliography}
\end{document}